\documentclass{article}
\usepackage{amsthm}
\usepackage{amsfonts}
\usepackage{amssymb}
\usepackage{amsmath}
\usepackage{esint} 
\usepackage{biblatex}
\usepackage{units} 
\newtheorem{thm}{Theorem}[section]
\newtheorem{corollary}[thm]{Corollary}
\newtheorem{lemma}[thm]{Lemma}
\newtheorem{proposition}[thm]{Proposition}

\theoremstyle{definition}
\newtheorem{defin}[thm]{Definition}
\newtheorem{remark}[thm]{Remark}
\newtheorem{example}[thm]{Example}

\newcommand*\Ca{\text{Cap}_1}
\newcommand*\interior[1]{\mathring{#1}}
\newcommand*\cls[1]{\overline{#1}}
\newcommand{\bb}[1]{\mathbb{#1}}
\newcommand{\ml}[1]{\mathcal{#1}}
\newcommand{\Emptyset}{\mathchar"001F}

\title{On the properties of the set where a generalized function of bounded variation takes infinite value}
\author{ \footnote{Mathematical Institute, University of Oxford, Oxford, UK. E-mail address: alessandro.cucinotta@maths.ox.ac.uk} Alessandro Cucinotta}

\begin{document}
\maketitle

\begin{abstract}
We study the properties of the set where a generalized function of bounded variation has infinite approximate limit, highlighting in this way the main geometric difference with functions of bounded variation. To this aim we prove a new result on strict approximation of sets of finite perimeter from the outside with open sets.
\end{abstract}


\section*{Introduction}
The main geometric difference between functions of bounded variation and generalized functions of bounded variation is that the latter may have infinite approximate limit on bigger sets. More precisely, the set of points where a function of bounded variation defined on $\Omega \subset \bb{R}^N$ has infinite approximate limit is negligible w.r.t. the $N-1$ dimensional Hausdorff measure, while it is easy to see that generalized functions of bounded variation may have infinite approximate limit on sets with Hausdorff dimension greater than $N-1$. In the paper we study this discrepancy.
\par
Let $L^0(\Omega)$ be the collection of Lebesgue-measurable functions that are finite valued almost everywhere. Given $u \in L^0(\Omega)$ we denote by $\tilde{u}(x)$ its approximate limit at $x$ whenever it exists and
 we say that a set $C \subset \Omega$ is $L^0$-polar if there exists a function $v \in L^0(\Omega)$ such that $C \subset \{x \in \Omega: |\tilde{v}(x)|=+ \infty \}$. Analogous definitions hold if $u$ is in the space $BV(\Omega)$ of functions of bounded variation or in the space $GBV(\Omega)$ of generalized functions of bounded variation, and we denote the collections of polar sets respectively with $\ml{P}_{L^0}$, $\ml{P}_{BV}$, and $\ml{P}_{GBV}$ (Definition \ref{D2}). \par 
Employing results from the general theory of approximate limits and functions of bounded variation we show that $\ml{P}_{L^0}$ is the class of sets that are negligible w.r.t. the Lebesgue measure $\lambda^N$ while $\ml{P}_{BV}$ is the class of sets that are negligible w.r.t. the Hausdorff measure $H^{N-1}$, and we prove that all inclusions in the following chain are strict:
\[
\ml{P}_{BV} \subset \ml{P}_{GBV} \subset \ml{P}_{L^0}.
\]
The class $\ml{P}_{GBV}$ cannot be characterized in terms of Hausdorff measures since one can show that $GBV$-polar sets may have any Hausdorff dimension smaller than the dimension of the ambient space. \par
To introduce our characterization of the class $\ml{P}_{GBV}$ we observe that the previous results may be reformulated saying that a set $C \subset \Omega$ belongs to $\ml{P}_{L^0}$ if and only if there exists a collection of open sets $\{U_k\}_{k \in \bb{N}}$ containing $C$ such that $\lambda^N(U_k) \to 0$, while it belongs to $\ml{P}_{BV}$ if and only if in addition to the previous conditions we also have that the perimeters $P(U_k,\Omega)$ of the sets $U_k$ tend to zero as $k$ goes to infinity.
We then prove that a subset $C \subset \Omega$ belongs to $\ml{P}_{GBV}$ if and only if 
there exists a collection $\{U_k\}_{k \in \bb{N}}$ of open sets of locally finite perimeter containing $C$ whose Lebesgue measure decreases to zero as $k$ goes to infinity. The proof of this fact cannot be obtained through the same tools used for $\ml{P}_{L^0}$ and $\ml{P}_{BV}$, and relies on a new approximation result for sets of finite perimeter which is interesting in itself.
\par 
More precisely we prove that a set of finite perimeter $C \subset \Omega$ can be approximated in a strict sense (i.e., in measure and with the perimeters of the approximating sets approaching the perimeter of $C$) by open sets containing the points of $\Omega$ where $C$ has density $1$. We also show that one cannot replace the set of points where $C$ has density $1$ with $C$ itself and that the approximating sets in general cannot have smooth boundary. This approximation theorem is not directly implied by existing ones (\cite[Theorem $3.42$]{AmbFuscPall}, \cite{Comi} and \cite{Schmidt}) and its proof relies on tools from capacity theory combined with a strong approximation result proved by Quentin de Gromard in \cite{grom}. \par 
The first section fixes the notation and contains preliminaries, in Section \ref{S2} we analyze the properties of polar sets, while the last section deals with the proof of the aforementioned approximation result.

\section{Notation and preliminaries}
Let $\Omega$ be an open subset of $\bb{R}^N$. Given $A \subset \bb{R}^N$ we denote its indicator function by $1_A$. We denote by $\omega_N$ the volume of the $N$-dimensional unit ball and given a Borel set $A \subset \Omega$ we will denote by $\lambda^N(A)$ its Lebesgue measure, while $H^{N-1}(A)$ will be its $N-1$ dimensional Hausdorff measure. Whenever we apply $\lambda^N$ or $H^{N-1}$ to a set we assume implicitly its $\lambda^N$- or $H^{N-1}$-measurability. We say that a sequence $\{\Omega_k\}_{k \in \bb{N}}$ of open subset of $\Omega$ is an exhaustion of $\Omega$ if for every $k$ we have $\cls{\Omega}_k \subset \Omega_{k+1}$ and the union of the sets $\Omega_k$ is the whole $\Omega$. It is well known that there exists an exhaustion of $\Omega$ made of smooth sets. \par
Given $x \in \Omega$ the upper and lower density of $A$ at $x$ are defined respectively as
\[
\limsup_{\rho \to 0^+} \frac{\lambda^N(A \cap B_\rho(x))}{\lambda^N(B_\rho(x))} 
\quad
\text{and}
\quad 
\liminf_{\rho \to 0^+} \frac{\lambda^N(A \cap B_\rho(x))}{\lambda^N(B_\rho(x))},
\]
while we say that $A$ has density $t$ at $x$ (and we write $\theta(A,x)=t$) if $t$ is both the upper and the lower density of $A$ at $x$.
The set of points in $\Omega$ where $A$ has density $t$ will be denoted by $A^{(t)}$, while the set of points where $A$ has strictly positive upper density will be denoted $A^+$. 
We will indicate with $L^0(\Omega)$ the collection of $\lambda^N$-a.e. real valued Lebesgue measurable functions and for any $u \in L^0(\Omega)$ its approximate upper limit at $x$ is
\[
u^+(x) := \inf \{t \in \mathbb{R} : \, \theta( \{ y \in \Omega : \, u(y)>t \}, x)=0 \},
\]
with the convention that $\inf({\Emptyset})= + \infty$. Similarly the approximate lower limit  of $u$ at $x$ is the value
\[
u^-(x) := \sup \{t \in \mathbb{R} : \, \theta( \{ y \in \Omega : \, u(y)<t \}, x)=0 \},
\]
with the convention that $\sup({\Emptyset})= - \infty$. If $u^+(x)=u^-(x)$ their common value is called the approximate limit of $u$ at $x$ and is denoted $\tilde{u}(x)$. It follows from the definition that for a continuous function $u$ we have $\tilde{u}(x)=u(x)$ for every $x \in \Omega$. It is also easy to see that if $f:\cls{\bb{R}} \to \cls{\bb{R}}$ is a continuous function and $u \in L^0(\Omega)$ has approximate limit $\tilde{u}(x)$ at $x$ then the approximate limit of $f \circ u$ at $x$ is $f(\tilde{u}(x))$. We will use the following theorem.

\begin{thm} \label{T3}
\cite[Theorem $2.9.13$]{GMT} If $u \in L^0(\Omega)$ then for $\lambda^N$-almost every $x \in \Omega$ the approximate limit of $u$ at $x$ exists and is finite.
\end{thm}
 
For the properties of the space $BV(\Omega)$ of functions of bounded variation we refer to \cite[Chapter $3$]{AmbFuscPall} and \cite[Chapter $5$]{EvansGari}. Given $u \in L^1_{loc}(\Omega)$ we denote its total variation (\cite[Definition $3.4$]{AmbFuscPall}) by $V(u,\Omega)$. If $u \in L^1(\Omega)$ then $u \in BV(\Omega)$ if and only if $V(u,\Omega)<+ \infty$. If $u \in BV(\Omega)$ then $Du$ will be its distributional derivative (which is a bounded Radon measure with values in $\bb{R}^N$) and $|Du|$ will be the variation of $Du$, which satisfies $|Du|(\Omega)=V(u,\Omega)$. It is known that the total variation is lower semicontinuous with respect to convergence in $L^1_{loc}(\Omega)$. Moreover, if $u \in BV(\Omega)$ and $f:\bb{R} \to \bb{R}$ is Lipschitz with $f(0)=0$, then $f \circ u \in BV(\Omega)$ and $|Df(u)|(\Omega)| \leq Lip(f)|Du|(\Omega)$.
We define the norm $\|u\|_{BV}:=\|u\|_{L^1}+|Du|(\Omega)$ and we refer to convergence in the $BV$ norm as strong convergence. \par
We will indicate the perimeter of $A$ in $\Omega$ (\cite[Definition $3.35$]{AmbFuscPall}) with $P(A, \Omega)$ and we say that $A$ has locally finite perimeter if its perimeter is finite in every precompact open subset of $\Omega$, in this case $\partial^*A$ will be its reduced boundary (\cite[Definition $3.54$]{AmbFuscPall}). If $ u \in BV(\Omega)$ the Coarea Formula (\cite[Theorem $3.40$]{AmbFuscPall}) implies that
\[
\int_{-\infty}^{+ \infty} P(\{u>t\},\Omega) \; dt =|Du|(\bb{R}^N).
\]
By the blow-up properties of the reduced boundary we have the following theorems.
\begin{thm} \cite[Theorem $3.59$]{AmbFuscPall}
Let $A$ be a set of finite perimeter in $\Omega$, then $P(A,\Omega)=H^{N-1}(\partial^*A)$.
\end{thm}

\begin{thm} \label{T2}
\cite[Theorem $3.61$]{AmbFuscPall}
Let $A$ be a set of finite perimeter in $\Omega$, then $\partial^*A \subset A^{(\nicefrac{1}{2})}$.
\end{thm}
We will use also the following fine properties of $BV$ functions.
\begin{thm} \label{T1} \cite[Theorems $2$-$3$, Section $5.9$]{EvansGari}
If $u \in BV(\Omega)$ then for  $H^{N-1}$-almost every $x \in \Omega$ we have that $u^+(x)$ and $u^-(x)$ are finite and
\[
\lim_{\rho \to 0}\fint_{B_\rho(x)}u(y) \; dy =\frac{u^+(x)+u^-(x)}{2},
\]
where the slashed integral denotes the mean value.
\end{thm}

The next result is due to Quentin de Gromard and will be crucial in Section \ref{S1}.
\begin{thm} \label{T7}
\cite[Theorem 3.1]{grom}
Let $B \subset \Omega$ be a set of finite perimeter in $\Omega$ and let $\epsilon >0$ be fixed. There exists a relatively closed set $L \subset \Omega$ such that the following properties hold:
\begin{enumerate}
\item $\lambda^N(B \Delta L) < \epsilon$;
\item $|D(1_B-1_L)|(\Omega) < \epsilon $;
\item $H^{N-1}((\partial L \cap \Omega) \setminus \partial^*L)< \epsilon$.
\end{enumerate}
\end{thm}
We now recall the definition of generalized functions of bounded variation.
\begin{defin} \label{D1}
A function $u: \Omega \rightarrow \bb{R}$ is a generalized function of bounded variation if for every $m \in \bb{R}_+$ the truncated function at the level $m$, i.e., $m \wedge u \vee -m$, belongs to $BV_{loc}(\Omega)$. We denote such truncation by $u^m$ and the space of these functions by $GBV(\Omega)$.
\end{defin}
If $u:\Omega \to \cls{\bb{R}}$ is real valued $\lambda^N$-a.e. and satisfies the truncation condition of Definition \ref{D1}, we say that $u \in GBV(\Omega)$, implicitly referring to its real valued representative.
Some properties of $GBV$ functions can be found in \cite[Section $4.5$]{AmbFuscPall}. Even if $GBV(\Omega)$ is not a vector space, the sum of positive $GBV$ functions is still in $GBV(\Omega)$. We will also use the following result.
\begin{proposition} \label{P13}
\cite[Theorem $4.34$]{AmbFuscPall}
Let $u \in GBV(\Omega)$, then for $\lambda^1$-a.e. $t \in \bb{R}$ the set $\{u>t\}$ has locally finite perimeter.
\end{proposition}
We now introduce the notation needed to apply the slicing techniques. Let $ \nu $ be a vector in $S^{N-1}$ and $C \subset \bb{R}^N$. We denote by $\pi_{\nu}$ the hyperplane of $\mathbb{R}^N$ orthogonal to $ \nu $ and by $C_{\nu}$ the orthogonal projection of $C$ on this hyperplane. 
For any $ y \in \pi_{\nu} $ the (possibly empty) set $\{ t \in \mathbb{R}: y+t \nu \in C \}$ is denoted by $ C_{\nu}^y $. \par
Given $u: C \rightarrow \cls{\bb{R}} $, for every $ y \in C_{\nu} $ such that $C_\nu^y \neq \emptyset$ the function $ u_{\nu}^y : C_{\nu}^y \rightarrow \cls{\bb{R}} $ is defined by $ u_{\nu}^y(t):=u(y+t \nu) $.
\begin{proposition} \label{P9} \cite[Proposition $4.35$]{AmbFuscPall}
Let $u \in GBV(\Omega)$ and $\nu \in S^{N-1}$, then for $H^{N-1}$-a.e. $y \in \Omega_\nu$ we have that $u_\nu^y \in GBV(\Omega_{\nu}^y)$ and
\[
(u^{\pm})_{\nu}^y(t)= (u_{\nu}^y)^{\pm}(t),
\]
for every $t \in \Omega_{\nu}^y $.
\end{proposition}

Finally we need some results about the $1$ capacity of a set. For these results we refer to \cite{FedZiem} and \cite[Section $4.7$]{EvansGari}. Using the notation of \cite{EvansGari} we say that a positive function $u$ belongs to $K^1$ if it is in $L^{\frac{N}{N-1}}(\bb{R}^N)$ and its distributional derivative is a vector valued function $\nabla u$ whose modulus is in $L^{1}(\bb{R}^N)$. 
\begin{defin}
The $1$-capacity of a set $E \subset \bb{R}^N$ is defined as the quantity
\[
\Ca(E):=
\]
\[
\inf \Big\{\int_{\bb{R}^N} |\nabla u (x)| \, dx: u \in K^1 , \; u \geq 1 \; \lambda^N \text{-a.e. on a neighborhood of } E \Big\},
\]
with the usual convention that $\inf \Emptyset:=+ \infty $.
\end{defin}
We will often refer to the $1$-capacity simply as the capacity. The next propositions summarize some of the main properties of the capacity.
\begin{proposition} \label{P15} \cite[Section $4.7$, Theorem $1$]{EvansGari}
The set function $\Ca (\cdot)$, defined on the power set of $\bb{R}^N$, is increasing and countably subadditive.
\end{proposition}

\begin{proposition} \label{P14} \cite[Section $5.6$, Theorem $3$]{EvansGari}
A set $E \subset \bb{R}^N$ satisfies $\Ca(E)=0$ if and only if $H^{N-1}(E)=0$.
\end{proposition}

The next proposition gives a characterization of the capacity that will be used in Proposition \ref{P7}.
\begin{proposition} \label{P6}
Given $E \subset \bb{R}^N$ the capacity of $E$ coincides with the following quantities:
\begin{enumerate}
\item[] $a(E)= \inf \Big\{P(B,\bb{R}^N):B \; \lambda^N \text{-measurable}, \, \lambda^N(B) < + \infty, \, E \subset \interior{B} \Big\}$,
\item[] $b(E) =\inf \Big\{|Du|(\bb{R}^N): u \in BV(\bb{R}^N) , \; u \geq 1 \; \lambda^N \text{-a.e. on a neigh. of } E \Big\}$.
\end{enumerate}
\begin{proof}
The equivalence of $\Ca(E)$ with $a(E)$ follows by \cite[Page $145$]{FedZiem}, so we only prove the equivalence of $a(E)$ and $b(E)$.
We only prove that $b(E) \geq a(E)$, since the other implication is trivial. To do this we prove that if $u \in BV(\bb{R}^N)$ and $u \geq 1$ $\lambda^N$-a.e. on a neighborhood of $E$, then there exists $B$ as in $a(E)$ such that $P(B,\bb{R}^N) \leq |Du|(\bb{R}^N)$.
To this aim we set $v \in BV(\bb{R}^N)$ as $v:=0 \vee u \wedge 1$ and we observe that $|Dv|(\bb{R}^N) \leq |Du|(\bb{R}^N)$. Moreover by the Coarea Formula we have that
\[
\int_0^1 P(\{v>t\},\bb{R}^N) \; dt =|Dv|(\bb{R}^N) \leq |Du|(\bb{R}^N),
\]
so that there exists $t_0 \in (0,1)$ such that $P(\{v>t_0\},\bb{R}^N) \leq |Du|(\bb{R}^N)$. \par
We now note that $\{v>t_0\} \supset \{u \geq 1\}$; since $u \geq 1$ at $\lambda^N$-a.e. every point of a neighborhood of $E$, there exists a $\lambda^N$-null set $A$ such that $\{v>t_0\} \cup A$ contains a neighborhood of $E$. Observe then that $P(\{v>t_0\} \cup A, \bb{R}^N)=P(\{v>t_0\}, \bb{R}^N)$ so that setting $B:=\{v>t_0\} \cup A$ we conclude.
\end{proof}
\end{proposition}

\begin{proposition} \label{P7}
There exists a dimensional constant $c >0$ such that for every $u \in BV(\bb{R}^N)$ and $\epsilon >0$ the following estimate is satisfied:
\[
\Ca\Big\{x: \exists \rho \in \bb{R}_+ : \; \fint_{B_\rho(x)}u(y) \; dy > \epsilon \Big\} \leq \frac{c}{\epsilon}|Du|(\bb{R}^N).
\]
\begin{proof}
The desired estimate is proved for every $u \in K^1$ in \cite[Lemma $1$, Section $4.8$]{EvansGari}. To prove our version one simply repeats the exact same argument of \cite{EvansGari} using the equivalence between $\Ca(E)$ and $b(E)$ proved in Proposition \ref{P6}.
\end{proof}
\end{proposition}

\section{Polar sets} \label{S2}
In this section we introduce the classes of polar sets $\ml{P}_{L^0}$, $\ml{P}_{BV}$, and $\ml{P}_{GBV}$. First we characterize $\ml{P}_{L^0}$ and $\ml{P}_{BV}$ respectively as the class of $\lambda^N$- and $H^{N-1}$-negligible sets and we observe that these conditions can both be expressed in terms of intersections of open sets (with a perimeter constraint in the $BV$ case). Then we use these characterizations to prove that the inclusions in the chain
$\ml{P}_{BV} \subset \ml{P}_{GBV} \subset \ml{P}_{L^0}$ are strict. To this aim we prove that any relatively closed $\lambda^N$-negligible set $C \subset \Omega$ belongs to $\ml{P}_{GBV}$ and that $GBV$-polar sets behave well w.r.t. one dimensional slicings. \par 
In the final part of the section, assuming a result which will be later proved in Section \ref{S1}, we characterize $\ml{P}_{GBV}$ in terms of intersections of open sets with locally finite perimeter (Theorem \ref{T4}), completing the picture on polar sets.

\begin{defin} \label{D2}
A set $C \subset \Omega$ is called $L^0$ polar (respectively $BV$ polar or $GBV$ polar) if there exists a function $u$ in $L^0(\Omega)$ (respectively in $BV(\Omega)$ or in $GBV(\Omega)$) such that $C \subset \{x \in \Omega : |\tilde{u}(x)|=+ \infty \}$. We denote the collection of these sets with $\ml{P}_{L^0}$ (respectively with $\ml{P}_{BV}$ or $\ml{P}_{GBV}$).
\end{defin}
Replacing $u$ with $|u|$ in the previous definition we obtain that a set $C \subset \Omega$ is ($L^0$, $BV$ or $GBV$) polar if and only if there exists a positive function $v$ (in $L^0(\Omega)$, $BV(\Omega)$ or $GBV(\Omega)$) such that $C \subset \{x \in \Omega : \tilde{v}(x)=+ \infty \}$.

\begin{proposition} \label{P1}
Let $C \subset \Omega$. The following conditions are equivalent:
\begin{enumerate}
\item[$(a)$] $C \in \ml{P}_{L^0}$;
\item[$(b)$] $\lambda^{N}(C)=0$;
\item[$(c)$] there exists a sequence of open sets $\{U_k\}_{k \in \bb{N}}$ containing $C$ such that $\lambda^N(U_k) \to 0$.
\end{enumerate}
\begin{proof}
If $C \in \ml{P}_{L^0}$ then $\lambda^N(C)=0$ by Theorem \ref{T3}. \par
If $\lambda^N(C)=0$ then there exists the desired sequence of open sets by the outer regularity of the Lebesgue measure. 
\par
Suppose now that we have a sequence as in $(c)$. Passing to a (not relabeled) subsequence $\{U_k\}_{i \in \bb{N}}$ we may suppose that  $\lambda^N(U_k) \leq k^{-2}$.
Now define $v:\Omega \rightarrow \bb{R} \cup \{+ \infty \}$ by
\[
v(x):= \sum_{k \in \bb{N}} 1_{U_k}(x)
\]
and observe that by monotone convergence theorem
\[
\int_{\Omega}|v(x)| \; dx = \sum_{k \in \bb{N}} \lambda^N(U_k) \leq 
\sum_{k \in \bb{N}}k^{-2} <+ \infty,
\]
so that $v$ is real valued $\lambda^N$-almost everywhere.
Moreover for every $k_0 \in \bb{N}$ we have that $v \geq k_0$ on $\cap_{k=1}^{k_0}U_k$, and since this is an open set containing $C$ we deduce that for every $x \in C$ we have $u^-(x) \geq k_0$. By the arbitrariness of $k_0$ we conclude.
\end{proof}
\end{proposition}

The next lemma is needed to characterize $\ml{P}_{BV}$ and will be used also in Section \ref{S1}. We denote by $S^N-1$ the unit sphere in $\bb{R}^N$ and by $\sigma_{N-1}$ its surface area.
\begin{lemma} \label{L1}
Let $Z \subset \Omega$ be such that $H^{N-1}(Z)<+ \infty$. There exists a dimensional constant $\tau$ such that for every $\epsilon >0$ there exists an open set  $V_{\epsilon} \supset Z$ with $\lambda^N(V_{\epsilon}) \leq \epsilon$ and $P(V_{\epsilon}, \Omega) \leq \tau(H^{N-1}(Z)+\epsilon)$.
\begin{proof}
To lighten the notation we set $c:=H^{N-1}(Z)$.
Let $\epsilon >0$ and $\delta >0$; by the definition of Hausdorff measure there exists a sequence of open balls $\{B_i\}_{i \in \bb{N}}$ each with radius $r_i$ less than $\delta$, such that their union contains $Z$ and
\begin{equation} \label{E3}
\frac{\omega_{N-1}}{2^{N-1}} \sum_{i \in \bb{N}}r_i^{N-1} \leq c+ \epsilon.
\end{equation}
Define now $
U_{\delta}:=\bigcup_{i \in \bb{N}}B_i \cap \Omega$
and note that this set is open and contains $Z$. Moreover, since each $r_i$ is less than $\delta$, taking \eqref{E3} into consideration, we get
\[
\lambda^N(U_\delta) \leq \sum_{i \in \bb{N}} \omega_N \, r_i^{N} \leq \delta \sum_{i \in \bb{N}} \omega_N \, r_i^{N-1}
 \leq \delta \, \frac{\omega_N 2^{N-1}}{\omega_{N-1}}(c+ \epsilon).
\]
Hence choosing $\delta$ small enough we have $\lambda^N(U_\delta)<\epsilon$. Reasoning similarly we obtain
\[
P(U_\delta, \Omega) \leq \sum_{i \in \bb{N}} P(B_i \cap \Omega, \Omega)
\leq \sigma_{N-1}\sum_{i \in \bb{N}}r_i^{N-1}
\leq \frac{2^{N-1}\sigma_{N-1}}{\omega_{N-1}}(c+\epsilon).
\]
In conclusion if $\delta $ is chosen small enough we can define $V_{\epsilon}:=U_{\delta}$.
\end{proof}
\end{lemma}

We are now ready to characterize the class $\ml{P}_{BV}$.
\begin{proposition} \label{P2}
Let $C \subset \Omega$. The following conditions are equivalent:
\begin{enumerate}
\item[$(a)$] $C \in \ml{P}_{BV}$;
\item[$(b)$] $H^{N-1}(C)=0$;
\item[$(c)$] there exists a sequence of open sets $\{U_k\}_{k \in \bb{N}}$ containing $C$ such that $\lambda^N(U_k) \to 0$ and $P(U_k,\Omega) \to 0$.
\end{enumerate}
\begin{proof}
If $C \in \ml{P}_{BV}$ then $H^{N-1}(C)=0$ by Theorem \ref{T1}, while if $H^{N-1}(C)=0$ there exists a sequence $\{U_k\}_{k \in \bb{N}}$ of open sets as in $c$ by Lemma \ref{L1}. \par 
Suppose now that there exists a sequence of open sets as in $(c)$. Passing to a (not relabeled) subsequence $\{U_k\}_{k \in \bb{N}}$ we may suppose that $\lambda^N(U_k) \leq k^{-2}$ and $P(U_k, \Omega) \leq k^{-2}$.
Reasoning as in Proposition \ref{P1} we define $v:\Omega \rightarrow \bb{R} \cup \{+ \infty \}$ by
\[
v(x):= \sum_{k \in \bb{N}} 1_{U_k}(x)
\]
and we have that $v \in L^1(\Omega)$ and $C \subset \{x:\tilde{v}(x)=+ \infty\}$. Moreover $v$ is the limit in $L^1(\Omega)$ of its partial sums, whose total variations are equibounded by the perimeter condition on $\{U_k\}_{k \in \bb{N}}$. By the lower semicontinuity of the variation we deduce that $v \in BV(\Omega)$.
\end{proof}
\end{proposition}

We now turn our attention to $\ml{P}_{GBV}$.
\begin{proposition} \label{P3}
Let $C \subset \Omega$ be a relatively closed set such that $\lambda^N(C)=0$, then $C \in \ml{P}_{GBV}$.
\begin{proof}
Consider $u: \Omega \rightarrow \cls{\bb{R}}$
defined by
\[
u(x):= \frac{1}{d(x,C)}.
\]
Since $C$ is a $\lambda^N$-null set then $u$ is real valued $\lambda^N$-almost everywhere. 
Fix $m >0$ and note that the truncated function $u^m$ satisfies
\[
u^m(x)= \frac{1}{d(x,C) \vee m^{-1}}.
\]
Being the reciprocal of a Lipschitz function strictly greater than $\frac{1}{m}$, the function $u^m$ is itself Lipschitz and belongs to $BV_{loc}(\Omega)$, implying that $u \in GBV(\Omega)$.
Moreover since $u$ is continuous we have that $\{y: |\tilde{u}(y)|=+ \infty \}=u^{-1}({+ \infty})=C $.
\end{proof}
\end{proposition}
The previous proposition implies that $GBV$-polar sets may have any Hausdorff dimension smaller than the dimension of the ambient space, so that the inclusion $\ml{P}_{BV} \subset \ml{P}_{GBV}$ is strict. The next proposition concerns $GBV$-polar sets in dimension $1$ and will be used to find a necessary condition for a set to be $GBV$-polar by means of a slicing argument.

\begin{proposition} \label{P4}
Let $\Omega \subset \bb{R}$ and let $C$ be a subset of $\Omega$. Then $C \in \ml{P}_{GBV} $ if and only if $\lambda^1(\cls{C} \cap \Omega)=0$.
\begin{proof}
By Proposition \ref{P3} we know that if $\lambda^1(\cls{C} \cap \Omega)=0$ then $C \in \ml{P}_{GBV}$. \par 
Viceversa suppose that $C \in \ml{P}_{GBV}$ and let $u \in GBV(\Omega)$ be a positive function such that $
C \subset \{y: \tilde{u}(y)=+ \infty \} $. Denote by $R_u$ the set of points where the approximate limit of $u$ exists and recall that $\lambda^1(\Omega \setminus R_u)=0$.
Fix $m \in \bb{N}$ and observe that $u^m:=m \wedge u $ has approximate limit $m$ at every point of $C$. We claim that its approximate limit is $m$ at every point of $\cls{C} \cap R_u$. \par 
To this aim consider the left continuous representative of $u^m$ (\cite[Theorem $3.28$]{AmbFuscPall}) and observe that it must have approximate limit $m$ at any point of $R_u$ (here $u^m$ has the approximate limit) that can be approximated from the left with points in $C$ (so that the approximate limit of $u^m$ must be greater than $m$).
Analogously the right continuous representative of $u^m$ will have approximate limit $m$ at any point of $R_u$ that can be approximated from the right with points in $C$. Since the approximate limit does not depend on
the representative we deduce that $u^m$ has approximate limit $m$ at every point of $\cls{C} \cap R_u$. Hence $u$ must have infinite approximate limit on every such point, implying by Theorem \ref{T3} that $\lambda^1(\cls{C} \cap \Omega)=\lambda^1(\cls{C} \cap R_u)=0$.
\end{proof}
\end{proposition}

\begin{proposition} \label{P5}
Let $\nu \in S^{N-1}$ and $C \in \ml{P}_{GBV}$. Then for $H^{N-1}$-a.e. $y \in C_\nu$ we have $\lambda^1(\cls{C_\nu^y} \cap \Omega_\nu^y)=0$.
\begin{proof}
By Proposition \ref{P9} for $H^{N-1}$-a.e. $y \in C_\nu$ we have $u_\nu^y \in GBV(\Omega_\nu^y)$ and 
\[
\{t \in \Omega_\nu^y : \widetilde{u_\nu^y}(t)=+ \infty\}=\{x \in \Omega: \tilde{u}(x)=+ \infty \}_\nu^y \supset C_\nu^y,
\]
so that Proposition \ref{P4} implies that
$\lambda^1(\cls{C_\nu^y} \cap \Omega_\nu^y)=0$.
\end{proof}
\end{proposition}

The previous proposition implies that the inclusion $\ml{P}_{GBV} \subset \ml{P}_{L^0}$ is strict. An example of a $L^0$-polar set which is not $GBV$-polar is the cartesian product of $\mathbb{Q}$ with any set $A \subset \bb{R}^{N-1}$ such that $\lambda^{N-1}(A)>0$, taking $\Omega:=\bb{R}^N$. \par 
The remaining part of the section is devoted to proving Theorem \ref{T4}, which characterizes $\ml{P}_{GBV}$ along the lines of conditions $(c)$ in Propositions \ref{P1} and \ref{P2}. We will use Proposition \ref{P11}, which follows by Theorem \ref{T5} (whose proof is postponed to Section \ref{S1}) by a standard localization argument.

\begin{proposition} \label{P11}
Let $A \subset \Omega$ be a set having locally finite perimeter in $\Omega$, then for every $\epsilon >0$ there exists $U \subset \Omega$, open set with locally finite perimeter in $\Omega$, such that $U \supset A^{(1)}$ and $\lambda^N(U \setminus A) \leq \epsilon$.
\end{proposition}

The following lemma is needed in view of Proposition \ref{P10}. 
\begin{lemma} \label{L5}
Let $C \in \ml{P}_{GBV}$ then there exists a positive function $u \in L^1(\Omega) \cap GBV(\Omega)$ such that $C \subset \{x \in \Omega:\tilde{u}(x)=+ \infty \}$.
\begin{proof}
Let $v \in GBV(\Omega)$ be a positive function such that $C \subset \{x \in \Omega :\tilde{v}(x)=+ \infty \}$.
Let then $\{\Omega_k\}_{k \in \bb{N}}$ be an exhaustion of $\Omega$ made of smooth sets and define $\Omega_0:= \Emptyset$. Now consider for every $k \geq 2$ the sets $B_k:=\Omega_k \setminus {\cls{\Omega}_{k-2}}$ and note that the product $v1_{B_k}$ is in $GBV(\Omega)$ for every $k$ and that $\lambda^N(B_k)<+ \infty$. As a consequence for every $k \in \bb{N}$ we can choose a sequence $\{t_i^k\}_{i \in \bb{N}} \subset \bb{R}_+$ such that $t_i \uparrow + \infty$ as $i$ goes to infinity and
\begin{equation} \label{E4}
P(\{v >t_i^k\} \cap B_k,\Omega)<+ \infty, \quad \lambda^N(\{v >t_i^k\} \cap B_k) < \frac{1}{2^i}.
\end{equation} 
We now define $v_k :\Omega \to \bb{R} \cup \{+ \infty \}$ by
\[
v_k(x):=\frac{1}{2^k}1_{B_k}(x)\sum_{i \in \bb{N}}1_{\{v>t_i^k\}}(x)
\]
and we note that because of \eqref{E4} we have $v_k \in GBV(\Omega)$ and
\begin{equation} \label{E5}
\int_\Omega |v_k(x)| \; dx \leq \frac{1}{2^k}.
\end{equation}
Moreover by construction $C \cap B_k \subset \{x \in \Omega:\widetilde{v_k}(x)=+ \infty \}$. Define now $u:\Omega \to \bb{R} \cup \{+ \infty \}$ by
\[
u:= \sum_{k \in \bb{N}}v_k
\]
and observe that this function is in $L^1(\Omega)$ by \eqref{E5} and that $u \in GBV(\Omega)$ as it is a locally finite sum of positive $GBV$ functions. Moreover if $x \in C$ then $x \in C \cap B_k$ for some $k \in \bb{N}$ so that $\tilde{u}(x)=+ \infty$ by construction.
\end{proof}
\end{lemma}

The next proposition provides a first characterization of $\ml{P}_{GBV}$. 
\begin{proposition} \label{P10}
Let $C \subset \Omega$, then $C \in \ml{P}_{GBV}$ if and only if there exists a sequence $\{A_k\}_{k \in \bb{N}}$ of  sets with locally finite perimeter in $\Omega$ such that $\lambda^{N}(A_k) \to 0$ and $\theta(A_k,x)=1$ for every $x \in C$ and for every $k \in \bb{N}$.

\begin{proof}
Suppose first that we have a collection $\{A_k\}_{k \in \bb{N}}$ as in the statement. It is not restrictive to assume $A_{k+1}\subset A_k$ for every $k \in \bb{N}$ (taking the intersection of the first $k$ sets). Then define $u: \Omega \rightarrow \bb{R} \cup{+ \infty}$ by
\[
u(x):=\sum_{k \in \bb{N}} 1_{A_k}(x)
\] 
and  note that $u$ is real valued $\lambda^N$-a.e. and that $ u \in GBV(\Omega)$ since its truncations are a finite sum of characteristic functions of sets having locally finite perimeter. 
If $x \in C$, for every $k \in \bb{N}$ we have
\[
\theta(\{u<k\},x) \leq \theta(\Omega \setminus A_k , x)=0,
\]
proving that $\tilde{u}(x)= + \infty$. Hence $C \in \ml{P}_{GBV}$.
\par 
Viceversa let $C \in \ml{P}_{GBV}$. By Lemma \ref{L5} there exists a positive function $u \in GBV(\Omega) \cap L^1(\Omega)$ such that $\tilde u(x)=+\infty$ for every $x \in C$. Note that this means that for every $t \in \bb{R}$ and for every $x \in C$ we have $\theta(\{u>t\},x)=1$. Since $u \in GBV(\Omega)$, by Proposition \ref{P13} we can find a sequence $\{t_k\}_{k \in \bb{N}}$ with $t_k\uparrow +\infty$ such that $\{u>t_k\}$ has locally finite perimeter in $\Omega$, and since $u \in L^1(\Omega)$ we also obtain that $\lambda^N(\{u>t_k\}) \downarrow 0$. Defining $A_k:=\{u>t_k\}$ we conclude.
\end{proof}
\end{proposition}

The next theorem, together with Propositions \ref{P1} and \ref{P2} completes the picture of polar sets.
\begin{thm} \label{T4}
Let $C \subset \Omega$, then $C \in \ml{P}_{GBV}$ if and only if there exists a sequence of open sets $\{U_k\}_{k \in \bb{N}}$ having locally finite perimeter in $\Omega$ and containing $C$ such that $\lambda^N (U_k) \to 0$.
\begin{proof}
If we have such a collection $\{U_k\}_{k \in \bb{N}}$ the statement follows by Proposition \ref{P10} since we have stronger hypotheses. \par
Viceversa suppose that $C \in \ml{P}_{GBV}$.  By Proposition \ref{P10} there exists a sequence $\{A_k\}_{k \in \bb{N}}$ of  sets with locally finite perimeter in $\Omega$ such that $\lambda^{N}(A_k) \to 0$ and $\theta(A_k,x)=1$ for every $x \in C$ and every $k \in \bb{N}$. By Proposition \ref{P11} for every $k \in \bb{N}$ there exists an open set $U_k$ with locally finite perimeter in $\Omega$ such that $U_k \supset A_k^{(1)}$ and $\lambda^N(U_k \setminus A_k) \leq \frac{1}{k}$.
This last fact, since $\lambda^N(A_k) \to 0$, implies that $\lambda^N(U_k) \to 0$. Finally for every $k \in \bb{N}$ we have that $U_k \supset A_k^{(1)} \supset C$, concluding the proof.
\end{proof}
\end{thm}

\section{Outer approximation of sets of finite perimeter with open sets} \label{S1}

The goal of the section is to prove that if $A \subset \Omega$ is a set of finite perimeter in $\Omega$, then for every $\epsilon >0$ there exists an open set $U$ such that $U \supset A^{(1)}$, $\lambda^N(U \setminus A) < \epsilon$, and $|P(U,\Omega)-P(A,\Omega)|<\epsilon$ (Theorem \ref{T5}). This fact then implies the result that we used in the previous section (Proposition \ref{P11}). To prove the aforementioned theorem we first need to show that strong convergence in $BV$ implies (up to passing to a subsequence) convergence $H^{N-1}$-almost everywhere (Theorem \ref{T6}). This fact seems to be known (for example it is considered in the general setting of metric spaces in \cite{Lahti}) but does not appear in the standard references about functions of bounded variation. Since it is an easy consequence of Proposition \ref{P7}, we give a complete proof. The next proposition is a preliminary version of Theorem \ref{T6}.

\begin{proposition} \label{P8}
Let $\{u_k\}_{k \in \bb{N}}$ be a sequence of positive functions of bounded variation such that $u_k \to 0$ strongly in $BV(\bb{R}^N)$. Then there is a (not relabeled) subsequence such that $u^+_k(x) \to 0$ for $H^{N-1}$-almost every $x \in \bb{R}$.
\begin{proof}
Passing to a subsequence we may suppose that for every $k \in \bb{N}$ we have $\|u_k\|_{BV} < 2^{-k}$ and with this extra hypothesis we fix $\epsilon >0$ and we prove that the set
\[
A:=\{x \in \bb{R}^N: \forall n \in \bb{N} \; \exists k \in \bb{N}:  k \geq n \text{ and } u_k^+(x)> 2\epsilon \}
\]
has $H^{N-1}$ zero measure. If we are able to do this then the statement follows by the arbitrariness of $\epsilon$. Since $H^{N-1}(A)=0$ if and only if $\Ca(A)=0$ by Proposition \ref{P14}, we prove the latter. To this aim observe that 
\[
A=\bigcap_{n \in \bb{N}} \bigcup_{k \geq n } \{x:u_k^+(x)> 2\epsilon \},
\]
so that the monotonicity and the subadditivity of the capacity (Proposition \ref{P15}) give that for every $n \in \bb{N}$ we have
\begin{equation} \label{E6}
\Ca(A) \leq \sum_{k \geq n}\Ca(\{x:u_k^+(x)> 2\epsilon \}).
\end{equation}
Taking into account that $u$ is positive together with Theorem \ref{T1} we get that
\[
\{x:u_k^+(x)> 2\epsilon \}
\subset
\{x:\frac{u_k^+(x)+u_k^-(x)}{2}> \epsilon \}
\]
\[
\subset
\{x: \exists \rho>0 : \fint_{B_\rho(x)}u_k(y) \; dy> \epsilon \},
\]
which implies by Proposition \ref{P7} that there exists a constant $C >0$ such that $\Ca(\{x:u_k^+(x)> 2\epsilon \}) \leq C|Du_k|(\bb{R}^n)$. Combining this with \eqref{E6} we deduce that for every $n \in \bb{N}$ we have
\[
\Ca(A) \leq C \sum_{k \geq n}\|u_k\|_{BV} \leq C \sum_{k \geq n} 2^{-k},
\]
and letting $n$ increase to infinity we get that $\Ca(A)=0$.
\end{proof}
\end{proposition}

\begin{thm} \label{T6}
Let $\{u_k\}_{k \in \bb{N}}$ be a sequence converging strongly in $BV(\bb{R}^N)$ to $u \in BV(\bb{R}^N)$. Then there exists a (not relabeled) subsequence such that $u_{k}^+(x) \rightarrow u^+(x)$
for $H^{N-1}$-almost every $x \in \bb{R}^N$.
\begin{proof}
We will prove that $|u^+_k-u^+| \leq |u_k-u|^+$
and applying Proposition \ref{P8} to (a subsequence of) $|u-u_k|$ we will then obtain that for $H^{N-1}$-almost every $x \in \Omega$
\[
\limsup_{k \rightarrow + \infty}|u^+_k(x)-u^+(x)| \leq \lim_{k \rightarrow + \infty} |u_k-u|^+(x)=0
\]
which implies the statement of the theorem. \par 
We now prove that $|u^+_k-u^+| \leq |u_k-u|^+$. By triangle inequality we have $u_k \leq u + |u_k-u|$ and passing to the approximate upper limit on both sides we get $u^+_k \leq (u + |u_k-u|)^+ \leq u^+ + |u_k-u|^+$, giving $u^+_k - u^+ \leq |u_k-u|^+$.
By symmetry we deduce that $|u^+_k-u^+| \leq |u_k-u|^+$.
\end{proof}
\end{thm}

The next corollary is obtained from the previous theorem by a standard localization argument.

\begin{corollary} \label{C1}
Let $\{u_k\}_{k \in \bb{N}} \subset BV_{loc}(\Omega)$ be a sequence converging strongly in $BV_{loc}(\Omega)$ to $u \in BV_{loc}(\Omega)$. There exists a (not relabeled) subsequence such that $u_k^+(x) \to u^+(x)$ for $H^{N-1}$-almost every $x \in \Omega$.
\end{corollary}

We recall that given $E \subset \Omega$ we denote by $E^{(1)}$ the set of points in $\Omega$ where $E$ has density $1$ and by $E^+$ the set of points where $E$ has positive upper density. By the Lebesgue Differentiation Theorem for $\lambda ^N$-almost every $x$ in $\Omega$ we have $1_E(x)=1_{E^+}(x)=1_{E^{(1)}}(x)$.
Moreover from the definition of upper density it easily follows that $
1^+_E(x)=1_{E^+}(x)$ and  that $E^{(1)}={^c}(({^c}E)^+)$, where ${^c}E:=\Omega \setminus E$. The next three lemmas provide the intermediate steps needed to prove Theorem \ref{T5}.

\begin{lemma} \label{L2}
Let $C \subset \Omega$ be a relatively closed set of finite perimeter in $\Omega$. Then
\[
C^+ \cup [(\Omega \cap \partial C) \setminus \partial^*C]=C.
\]
\begin{proof}
Note that by Theorem \ref{T2} we have $C^+ \supset \interior{C} \cup \partial^*C $, so that
\[
C^+ \cup [(\Omega \cap \partial C) \setminus \partial^*C]
\supset \interior{C} \cup (\Omega \cap \partial C) \supset C.
\]
Viceversa since $C$ is relatively closed we have $C^+ \cup (\Omega \cap \partial C) \subset C $, so that also $C^+ \cup [(\Omega \cap \partial C) \setminus \partial^*C] \subset C$.
\end{proof}
\end{lemma}

\begin{lemma} \label{L4}
Let $B \subset \Omega$ be a set of finite perimeter in $\Omega$. For every $\epsilon >0$ there exist a sequence $\{C_k\}_{k \in \bb{N}}$ of relatively closed sets of finite perimeter in $\Omega$ and a Borel set $Z \subset \Omega$ with the following properties:
\begin{enumerate}
\item[$(a)$] $|D(1_B-1_{C_k})|(\Omega) \to 0$ and $\lambda^N(B \Delta C_k) \to 0$;
\item[$(b)$] for every $x \notin Z$ we have that $1_{C_k}(x) \to 1_{B^+}(x)$;
\item[$(c)$] $H^{N-1}(Z)< \epsilon$.
\end{enumerate}
\begin{proof}
By Theorem \ref{T7} there exists a sequence $\{C_k \}_{k \in \bb{N}}$ of relatively closed sets with finite perimeter in $\Omega$ such that $\lambda^N(C_k \Delta B) \to 0$, $|D(1_{C_k}-1_{B})|(\Omega) \to 0$ and
\begin{equation} \label{E8}
H^{N-1}((\partial C_k \cap \Omega) \setminus \partial^* C_k) < \frac{\epsilon}{2^k}.
\end{equation}
By Corollary \ref{C1} there exists a $H^{N-1}$-null set $\tilde{Z} \subset \Omega$ and a (not relabeled) subsequence of $\{C_k \}_{k \in \bb{N}}$ such that for every $x \notin \tilde{Z}$ we have $1_{C_k}^+ (x) \rightarrow 1_B^+ (x)$ and, taking into account that $1^+_{A}=1_{A^+}$, we deduce that for every such $x$ we have 
\begin{equation} \label{E1}
1_{C_k^+} (x) \rightarrow 1_{B^+} (x).
\end{equation}
Define now 
\[
Z:=\tilde{Z} \cup  \bigcup_{k \in \bb{N}} [(\partial C_k \cap \Omega) \setminus \partial^*C_k]
\]
and note that by \eqref{E8} we have that $H^{N-1}(Z) < \epsilon$.
Since every $C_k$ is relatively closed in $\Omega$, by Lemma \ref{L2} we have that $C_k^+ \setminus Z=C_k \setminus Z$, which together with \eqref{E1} implies that for every $x \notin Z$ we have $1_{C_k}(x) \rightarrow 1_{B^+}(x)$.
\end{proof}
\end{lemma}

\begin{lemma} \label{L3}
Let $A \subset \Omega$ be a set of finite perimeter in $\Omega$. Then for every $\epsilon >0$ there exist $Z \subset \Omega$ such that $H^{N-1}(Z)< \epsilon$ and an open set $U$ such that $U \cup Z \supset A^{(1)}$, $\lambda^N(U \setminus A) < \epsilon$, and $P(U,\Omega)<P(A,\Omega)+\epsilon$.
\begin{proof}
Let $B:= \Omega \setminus A$ and note that $P(B,\Omega)=P(A,\Omega)< + \infty $. By Lemma \ref{L4} there exists sequence $\{C_k \}_{k \in \bb{N}}$ of relatively closed sets with finite perimeter in $\Omega$ and a Borel set $Z \subset \Omega$ such that: 
\begin{enumerate}
\item[$(a)$] $|D(1_B-1_{C_k})|(\Omega) < \epsilon 2^{-k}$ and $\lambda^N(B \Delta C_k) < \epsilon 2^{-k}$;
\item[$(b)$] For every $x \notin Z$ we have $1_{C_k}(x) \to 1_{B^+}(x)$;
\item[$(c)$] $H^{N-1}(Z)< \epsilon$.
\end{enumerate}
We define 
\[
C:=\bigcap_{k \in \bb{N}} C_k
\]
and $U:=\Omega \setminus C$, and we claim that $U$ has the required properties.
By $(b)$ we get that $
 B^+ \setminus Z \supset C \setminus Z $ and passing to the complement in this inclusion, keeping in mind that $A^{(1)}={^c}(B^+)$, we get $A^{(1)} \subset Z \cup U$.
Moreover
\[
\lambda^N(U \setminus A) \leq \sum_{k \in \bb{N}} \lambda^N((\Omega \setminus C_k) \setminus A) \leq 
\sum_{k \in \bb{N}} \lambda^N((\Omega \setminus A) \Delta C_k) =
\sum_{k \in \bb{N}} \lambda^N(C_k \Delta B),
\]
and the last term of the chain is less than $\epsilon$ by $(a)$. \par
To conclude the proof we only need to show that $P(U,\Omega) < P(A,\Omega)+ \epsilon$. To achieve this we prove that the function $1_U-1_A$ is in $BV_{loc}(\Omega)$ and satisfies $|D(1_{U}-1_A)|(\Omega) < \epsilon$,
so that we obtain 
\[
P(U,\Omega)=|D1_{U}|(\Omega) \leq |D(1_{U}-1_A)|(\Omega)+|D1_A|(\Omega)< P(A,\Omega)+ \epsilon.
\]
Observe that $1_U-1_A=(1-1_C)-(1-1_B)=1_B-1_C$ so that it is sufficient to show that $1_B-1_C \in BV_{loc}(\Omega)$ and satisfies $|D(1_B-1_C)|(\Omega) < \epsilon$. \par 
To this aim let $f: \bb{R} \rightarrow \bb{R}$ be the truncation function at the levels zero and one, i.e. $f(t):=0 \vee t \wedge 1$, and we claim that the following identity holds $\lambda^N$-almost everywhere in $\Omega$:
\begin{equation} \label{E2}
1_{B}-1_{C}=f\Big(\sum_{k \in \bb{N}} (1_B-1_{C_k})\Big).
\end{equation}
First of all observe that for every $x \notin Z$ the term $1_{B^+}(x)-1_{C_k}(x)$ is eventually zero by the condition $(b)$, so that since $\lambda^N(Z)=0$ and $\lambda^N(B \Delta B_+)=0$ we have that the series in the r.h.s. of \eqref{E2} is well defined $\lambda^N$-almost everywhere.
Another consequence of $(b)$ is that $\lambda^N$-almost every point of $C$ belongs to $B^+$ so that $\lambda^N$-almost every such point is in $B$. Taking this into account \eqref{E2} follows.
As a consequence the function $1_B-1_C$ is the limit in $L^1_{loc}(\Omega)$ of the sequence $\{f_j\}_{j \in \bb{N}}$ given by
\[
f_j:=f\Big(\sum_{k = 1}^{j} (1_B-1_{C_k})\Big),
\]
and since $f$ is Lipschitz with Lipschitz constant one, each function $f_j$ belongs to $BV_{loc}(\Omega)$ and satisfies
\[
|D f_j|(\Omega)
\leq \ \sum_{k \in \bb{N}} |D(1_B-1_{C_k})|(\Omega)< \epsilon.
\]
This implies that also $1_B-1_C$ belongs to $BV_{loc}(\Omega)$ and satisfies $|D(1_B-1_C)|(\Omega)< \epsilon$, concluding the proof.
\end{proof}
\end{lemma}

\begin{thm} \label{T5}
Let $A \subset \Omega$ be a set of finite perimeter in $\Omega$. Then for every $\epsilon >0$ there exists an open set $U$ such that $U \supset A^{(1)}$, $\lambda^N(U \setminus A) < \epsilon$, and $|P(U,\Omega)-P(A,\Omega)|<\epsilon$.
\begin{proof}
By Lemma \ref{L3} for every $n \in \bb{N}$ we find an open set $V_n$ and a set $Z_n$ such that $V_n \cup Z_n \supset A^{(1)}$, $\lambda^N(V_n \setminus A)< \frac{1}{n}$, $P(V_n,\Omega) < P(A,\Omega)+ \frac{1}{n}$ and $H^{N-1}(Z)< \frac{1}{\tau n}$, where $\tau$ is the dimensional constant of Lemma \ref{L1}. \par
By the aforementioned lemma we then find an open set $W_n \supset Z_n$ such that $ \lambda^N(W_n)< \frac{1}{n}$ and $P(W_n,\Omega)< \frac{1}{n}$,  and we define $U_n:=V_n \cup W_n$.
In this way we have a sequence $\{U_n\}_{n \in \bb{N}}$ of open sets containing $A^{(1)}$, converging to $A$ in measure and such that $P(U_n,\Omega) \leq P(V_n,\Omega)+P(W_n,\Omega) \leq P(A,\Omega)+ \frac{2}{n}$. Then by the lower semicontinuity of the perimeter w.r.t. convergence in measure it is sufficient to define $U:=U_n$ for $n$ sufficiently large.
\end{proof}
\end{thm}

\begin{remark}
The previous theorem makes sense even if the set $A$ is defined modulo $\lambda^N$-negligible sets, as all the terms involved are invariant under modifications on $\lambda^N$-null sets.
\end{remark}

\begin{remark}
Theorem \ref{T5} fails if we replace the condition that $U \supset A^{(1)}$ with $U \supset A$. We show this by considering the case when $\lambda^N(A)=0$. Indeed, if we could find a sequence of open sets $\{U_n\}$ containing $A$ with Lebesgue measure and perimeter decreasing to zero, we would obtain by Proposition \ref{P2} that $A \in \ml{P}_{BV}$, while the same proposition shows that this is in general not possible.
\end{remark}

The next example shows that we cannot require the smoothness of the approximating set in Theorem \ref{T5} when $N \geq 2$.
\begin{example}
Let $\Omega:=\bb{R}^N$ and let $A \subset \bb{R}^N$ be an open dense set with finite perimeter and finite Lebesgue measure (for example a countable union of balls with dense centers and sufficiently small radii). Suppose by contradiction that there exists a smooth open set $U$ such that $U \supset A^{(1)} \supset A$ and $\lambda^N(U)< + \infty$. If we had such a set $U$ we would then get that $\cls{U}=\bb{R}^N$, so that $\lambda^N(\partial U)=\lambda^N(\cls{U})-\lambda^N(U)=+ \infty$. On the other hand the smoothness of $U$ implies that $\lambda^N(\partial U)=0$, giving the desired contradiction.
\end{example}

\section*{Acknowledgement}
I would like to thank Prof. Gianni Dal Maso for his valuable advice. Most results of this paper were first obtained in the context of my master's thesis for the joint program of the University of Trieste and SISSA.

\printbibliography
\end{document}